\documentclass[a4paper]{amsart}
\usepackage{amsmath,amsthm,amssymb}
\usepackage[dvipdfmx]{graphicx}
\usepackage{mathrsfs}
\usepackage[dvipdfmx, usenames]{color}

\newtheorem{thm}{Theorem}[section]
\newtheorem{lem}[thm]{Lemma}
\newtheorem{cor}[thm]{Corollary}
\newtheorem{prop}[thm]{Proposition}

\theoremstyle{definition}

\newtheorem{defn}[thm]{Definition}

\theoremstyle{remark}

\newtheorem{rem}[thm]{Remark}        

\numberwithin{equation}{section}

\def\XXint#1#2#3{{\setbox0=\hbox{$#1{#2#3}{\int}$}
\vcenter{\hbox{$#2#3$}}\kern-.5\wd0}}

\newcommand{\Pro}{\mathcal{P}}
\newcommand{\R}{\mathbb{R}}

\newcommand{\Cpl}{\mathfrak{Cpl}}
\newcommand{\di}{\mathrm{Diam}}
\newcommand{\supp}{\mathrm{supp}\,}

\newcommand{\RE}{\mathcal{R}}
\newcommand{\LIP}{{\bf \mathrm{LIP}}}
\newcommand{\Ch}{\mathsf{Ch}}

\newcommand{\vol}{\mathrm{vol}}

\newcommand{\calE}{\mathcal{E}}

\newcommand{\calH}{\mathcal{H}}

\newcommand{\calL}{\mathcal{L}}
\newcommand{\calM}{\mathcal{M}}

\newcommand{\calR}{\mathcal{R}}

\newcommand{\calW}{\mathcal{W}}

\newcommand{\bbN}{\mathbb{N}}

\newcommand{\bbS}{\mathbb{S}}

\usepackage[abbrev,nobysame]{amsrefs}

\makeatletter
\def\@makefnmark{%
\leavevmode
\raise.9ex\hbox{\check@mathfonts
\fontsize\sf@size\z@\normalfont%
\@thefnmark}%
}
\makeatother

\title[Bishop-type inequality on $RCD$ spaces]{A Bishop type inequality on metric measure spaces with Ricci curvature bounded below}
\author{Yu Kitabeppu}
\address[Yu Kitabeppu]{Kyoto University}
\email{y.kitabeppu@gmail.com}
\begin{document}
\maketitle
 \begin{abstract}
  We define a Bishop-type inequality on metric measure spaces with Riemannian curvature-dimension condition. The main result in this short article is that any $RCD$ spaces with the Bishop-type inequalities possess only one regular set in not only the measure theoretical sense but also the set theoretical one. As a corollary, the Hausdorff dimension of such $RCD^*(K,N)$ spaces are exactly $N$. We also prove that every tangent cone at any point on such $RCD$ spaces is a metric cone. 
 \end{abstract}
 %
 %
\section{Introduction}
Ricci limit spaces are the Gromov-Hausdorff limits of sequences of Riemannian manifolds with Ricci curvature bounded from below. They are roughly classified into two classes, \emph{collapsing} and \emph{non-collapsing}. For simplicity, we only consider closed Riemannian manifolds in this section. Let $\{(M_n,g_n)\}_{n\in\bbN}$ be a sequence of $N$-dimensional closed Riemannian manifolds with Ricci curvature bounded from below and uniformly bounded diameters. By the Gromov compactness theorem, we are able to find a convergent subsequence $\{(M_{n_k},g_{n_k})\}_k$. We denote the Riemannian volume measure of $(M_n,g_n)$ by $\vol_n$. Let $(Y,\nu)$ be the limit space of such subsequence $\{(M_{n_k},g_{n_k})\}$, more precisely, of $\{(M_{n_k},g_{n_k},\underline{\vol}_{n_k})\}_k$ where $\underline{\vol}$ is the normalized volume measure $\vol/\vol(M)$. Then we call $(Y,\nu)$ a \emph{non-collapsing} Ricci limit space if 
\begin{align}
 \lim_{n\rightarrow\infty}\vol_n(M_n)>0,\notag
\end{align}
and we call $(Y,\nu)$ a \emph{collapsing} Ricci limit space if $\lim_{n\rightarrow\infty}\vol_n(M_n)=0$. On the series of papers of Colding and Cheeger-Colding \cite{Covolume,CC1,CC2,CC3}, it is proven that there are many nice properties on non-collapsing Ricci limit spaces.
 \begin{thm}[\cite{CC1,CC3,Covolume}]\label{thm:nonricci}
 Let $(Y,\nu)$ be a non-collapsing Ricci limit space of a sequence of $N$-dimensional Riemannian manifolds. Then  
\begin{itemize}
 \item $\dim_{\calH}Y=N$,
 \item $\nu$-almost every point has the unique tangent cone that is isometric to $\R^N$, 
 \item $\nu$ is absolutely continuous to $\calH^N$ on some nice regular set of full measure and vice versa,
 \item all the tangent cones are metric cones, 
 \item $(Y,\nu)$ is rectifiable in the sense of Cheeger-Colding $($see Definition $5.3$ in \cite{CC3}$)$. 
\end{itemize}
\end{thm}
Compare with the collapsing and non-collapsing Ricci limit spaces.  
\begin{thm}[\cite{CC1,CC3,CNholder}]\label{thm:colricci}
 Let $(Y,\nu)$ be a collapsing Ricci limit space of a sequence of $N$-dimensional Riemannian manifolds. Then 
 \begin{itemize}
  \item $\dim_{\calH}Y\leq N-1$, 
  \item there exists an integer $k$ with $1\leq k< N$ such that $\nu$-almost every points has the unique tangent cone that is isometric to $\R^k$,
  \item $\nu$ is absolutely continuous to $\calH^k$ on some nice regular set of full measure and vice versa,
  \item $(Y,\nu)$ is rectifiable in the sense of Cheeger-Colding. 
 \end{itemize}
\end{thm}
It is known that there is an example of collapsing Ricci limit space on which there exists a point whose tangent cones are not metric cones (Example 8.95 in \cite{CC1}). 
\par On the other hand, without using the Gromov-Hausdorff approximation by Riemannian manifolds, a concept of a metric measure space with lower Ricci curvature bound is also defined. Lott-Villani, Sturm have introduced the \emph{curvature-dimension condition} $CD(K,N)$ independently for the case $K\in\R$, $N=\infty$ and $K=0$, $1<N<\infty$ by \cite{LV}, for any case by \cite{Stmms1,Stmms2}. They prove that any Ricci limit spaces with finite total mass satisfy $CD$ conditions. However Finsler manifolds with $N$-Ricci curvature bounded below from $K$ satisfies $CD(K,N)$ condition though generic such manifolds are not Ricci limit spaces (see \cite{OhFinsler} for the definition of $N$-Ricci curvature and the relation with $CD$ condition, and see \cite{CC1} for the relation of Ricci limit spaces and Finsler manifolds). In order to rule out Finsler manifolds, Ambrosio-Gigli-Savar\'e introduce the \emph{Riemannian curvature-dimension condition} $RCD(K,\infty)$ for $K\in\R$ on metric measure space with finite mass \cite{AGSRiem} and with $\sigma$-finite measure \cite{AGMR}. Later Erbar-Kuwada-Sturm \cite{EKS} and Ambrosio-Mondino-Savar\'e \cite{AMS} define $RCD^*(K,N)$ condition for finite $N\in(1,\infty)$ independently, where $~^*$ means the \emph{reduced curvature-dimension condition} that is defined by Bacher-Sturm \cite{BS} for the tensorization property and the local-to-global property. Recently, Mondino-Naber prove the nice local structures on finite dimensional $RCD$ spaces \cite{MN}. 
\begin{thm}[\cite{MN,Stmms2}]
 Let $(X,d,m)$ be an $RCD^*(K,N)$ space for $K\in\R$ and $N\in(1,\infty)$. Assume $\supp m=X$.Then 
 \begin{itemize}
  \item $\dim_{\calH}(X,d)\in(0,N]$,
  \item for $m$-almost every point $x\in X$, there exists an integer $k=k_x$ with $1\leq k\leq N$ such that the tangent cones are unique and isomorphic to $\R^k$, 
  \item $(X,d,m)$ is rectifiable in the sense weaker than Cheeger-Colding's one.  
 \end{itemize}
\end{thm} 
We denote by $\calR_k$ the set of all points whose tangent cones are isomorphic to $k$-dimensional Euclidean space, and by $\calR=\cup_k\calR_k$ (see Definition \ref{def:regularsets}). 
\smallskip
\par Since $RCD$ spaces are not defined as limits of sequences of Riemannian manifolds, neither \emph{collapsing} nor \emph{non-collapsing} makes sense. In this article, we prove the $RCD$ spaces with Bishop-like inequalities behave like non-collapsing spaces. More precisely, our main result is as follows. 
\begin{thm}[Theorem \ref{thm:noncollapsing}, Corollary \ref{cor:ulbound}, Corollary \ref{cor:Hdim}, and Theorem \ref{thm:metriccone}]
 Let $(X,d,m)$ be an $RCD^*(K,N)$ space with the generalized Bishop inequality $BI(K,N)$. Assume $\supp m=X$. Then 
 \begin{enumerate}
  \item $\calR=\calR_N$, 
  \item on $\calR_N$, the measure $m$ is absolutely continuous to $\calH^N$ and vice versa,  
  \item $\dim_{\calH}(X,d)=N$,
  \item All tangent cones are metric cones. 
 \end{enumerate}
\end{thm}
\begin{rem}
 The family of $RCD^*(K,N)$ spaces with the generalized Bishop inequality $BI(K,N)$ includes the family of non-collapsing Ricci limit spaces of appropriate curvature and dimension bounds. However these two classes do not coincide. See Remark \ref{rem:boundaryTan}. 
\end{rem}
 %
 %
\section{Preliminaries}\label{sec:pre}
A triplet $(X,d,m)$ consisting of a complete separable metric space $(X,d)$ and a locally finite $\sigma$-compact positive Borel measure $m$ on $X$ is called a \emph{metric measure space}. Two metric measure spaces $(X,d,m)$ and $(Y,r,\nu)$ are \emph{isomorphic} if there exists an isometry $f:\supp m\rightarrow\supp \nu$ with $f_*m=\nu$. 
A continuous curve $\gamma:[0,1]\rightarrow X$ is \emph{absolutely continuous} if there exists an $L^1(0,1)$ function $g$ such that 
\begin{align}
  d(\gamma_s,\gamma_t)\leq \int_s^tg(r)\,dr,\qquad\text{for any }s\leq t\in[0,1].\notag
\end{align}
For a continuous curve $\gamma:[0,1]\rightarrow X$, the metric derivative $\vert\dot\gamma\vert$ is defined by 
\begin{align}
 \vert\dot\gamma_t\vert:=\lim_{s\rightarrow t}\frac{d(\gamma_t,\gamma_s)}{\vert s-t\vert}\notag
\end{align}
as long as the right-hand side makes sense.  
It is known that every absolutely continuous curve has the metric derivative for almost every point \cite{AGSbook}. We call an absolutely continuous curve $\gamma:[0,1]\rightarrow X$ a \emph{geodesic} if $\vert\dot\gamma_t\vert=d(\gamma_0,\gamma_1)$ for almost every $t\in [0,1]$. A metric space $(X,d)$ is called a \emph{geodesic space} if for any two points, there exists a geodesic connecting them. 

We denote the set of all Lipschitz functions on $X$ by $\LIP(X)$. For $f\in\LIP(X)$, the local Lipschitz constant at $x$, $\vert \nabla f\vert(x)$, is defined as
\begin{align} 
 \vert \nabla f\vert(x):=\limsup_{y\rightarrow x}\frac{\vert f(x)-f(y)\vert}{d(x,y)}\notag
\end{align} 
if $x$ is not isolated, otherwise $\vert \nabla f\vert(x)=\infty$. For $f\in L^2(X,m)$, we define the \emph{Cheeger energy} $\Ch(f)$ as 
\begin{align}
 \Ch(f):=\frac{1}{2}\inf\left\{\liminf_{n\rightarrow \infty}\int_X\vert \nabla f_n\vert^2\,dm\;;\;f_n\in\LIP(X),\;f_n\rightarrow f\text{ in }L^2(X,m)\right\}.\notag
\end{align}
Set $D(\Ch):=\{f\in L^2(X,m)\;;\;\Ch(f)<\infty\}$. We define the \emph{Sobolev space} $W^{1,2}(X,d,m):=L^2(X,m)\cap D(\Ch)$ equipped with the norm $\Vert f\Vert_{1,2}^2:=\Vert f\Vert_2^2+2\Ch(f)$. It is known that $W^{1,2}$ is a Banach space. We say that $(X,d,m)$ is \emph{infinitesimally Hilbertian} if $W^{1,2}$ is a Hilbert space. 

We denote the set of all Borel probability measures on $X$ by $\Pro(X)$. We define $\Pro_2(X)$ as the set of all Borel probability measures with finite second moment, that is, $\mu\in\Pro_2(X)$ if and only if $\mu\in\Pro(X)$ and there exists a point $o\in X$ such that $\int_Xd(x,o)^2\,d\mu(x)<\infty$. We call a measure $q\in\Pro(X\times X)$ a \emph{coupling} between $\mu$ and $\nu$ if $(p_1)_*q=\mu$ and $(p_2)_*q=\nu$, where $p_i:X\times X\rightarrow X$ are natural projections for $i=1,2$. For two probability measures $\mu_0,\mu_1\in\Pro_2(X)$, we define the \emph{$L^2$-Wasserstein distance} between $\mu_0$ and $\mu_1$ as 
\begin{align}
 W_2(\mu_0,\mu_1):=\inf\left\{\int_{X\times X}d(x,y)^2\,dq(x,y)\;;\;q\in\Cpl(\mu_0,\mu_1)\right\}^{1/2},\notag
\end{align}
where $\Cpl(\mu_0,\mu_1)$ is the set of all couplings of $\mu_0$ and $\mu_1$. The pair $(\Pro_2(X),W_2)$ is called the \emph{$L^2$-Wasserstein space}, which is a complete separable geodesic metric space if so is $(X,d)$. We explain how geodesics in $X$ relates to those in $\Pro_2(X)$. We denote the space of all geodesics in $X$ by $Geo(X)$, equipped with the sup distance. Define the \emph{evaluation map} $e_t:Geo(X)\rightarrow X$ for $t\in [0,1]$ by $e_t(\gamma)=\gamma_t$. Let $(\mu_t)_t\in Geo(\Pro_2(X))$ be a geodesic connecting $\mu_0,\mu_1$ in $\Pro_2(X)$. Then there exists a probability measure $\pi\in\Pro(Geo(X))$ such that $(e_t)_*\pi=\mu_t$, by which we say that the geodesic $(\mu_t)_t$ can be lifted to $\pi$. 
\subsection{The curvature-dimension condition}
 For given $K\in\R$ and $N\in(1,\infty)$, we define the \emph{distortion coefficients}, $\sigma^{(t)}_{K,N}$ for $t\in[0,1]$, by 
 \begin{align}
  \sigma^{(t)}_{K,N}(\theta):=\begin{cases}
                                            \infty&\text{if }K\theta^2\geq N\pi^2,\\
                                            \frac{\sin(t\theta\sqrt{K/N})}{\sin(\theta\sqrt{K/N})}&\text{if }0<K\theta^2<N\pi^2,\\
                                            t&\text{if }K\theta^2=0,\\
                                            \frac{\sinh(t\theta\sqrt{-K/N})}{\sinh(\theta\sqrt{-K/N})}&\text{if }K\theta^2<0.
                                           \end{cases}\notag
 \end{align}
 Let $(Y,d)$ be a geodesic metric space and $f:Y\rightarrow\R\cup\{\pm\infty\}$ a function on $Y$. 
 \begin{defn}[\cite{EKS}]
  A function $f:Y\rightarrow \R\cup\{\pm\infty\}$ is said to be \emph{$(K,N)$-convex} for $K\in\R$ and $N\in(1,\infty)$ if for any two points $y_0,y_1\in Y$, there exists a geodesic $(y_t)_t$ connecting them such that 
  \begin{align}
   &\exp\left(-\frac{1}{N}f(y_t)\right)\notag\\
   &\geq \sigma^{(1-t)}_{K,N}(d(y_0,y_1))\exp\left(-\frac{1}{N}f(y_0)\right)+\sigma^{(t)}_{K,N}(d(y_0,y_1))\exp\left(-\frac{1}{N}f(y_1)\right)\notag
  \end{align} 
  holds for any $t\in[0,1]$. 
 \end{defn} 
 Let $(X,d,m)$ be a geodesic metric measure space. Consider $\mu=\rho m\ll m$ a probability measure that is absolutely continuous with respect to $m$ and its Radon-Nikodym derivative being $\rho$. We define the relative entropy functional $\mathrm{Ent}_m$ by 
 \begin{align}
  \mathrm{Ent}_m(\mu):=\int_{\{\rho>0\}}\rho\log\rho\,dm,\notag
 \end{align}
 whenever $(\rho\log\rho)_+$ is integrable, otherwise we define $\mathrm{Ent}_m(\mu)=\infty$. 
%
 \begin{defn}[\cite{EKS}, cf. \cite{AMS}]
  Let $(X,d,m)$ be a geodesic metric measure space. We say that $(X,d,m)$ satisfies the \emph{entropic curvature-dimension condition} $CD^e(K,N)$ for $K\in\R$ and $N\in(1,\infty)$ if the relative entropy functional $\mathrm{Ent}_m$ is $(K,N)$-convex. Moreover if $(X,d,m)$ is infinitesimally Hilbertian, $(X,d,m)$ is called an $RCD^*(K,N)$ space. 
 \end{defn}
 Under the infinitesimal Hilbertianity condition, $CD^e(K,N)$ is equivalent to $CD^*(K,N)$. 
 \subsection{Tangent cones and regular sets on $RCD$ spaces}
 Let $(X,d,m)$ be a metric measure space. Take a point $x_0\in \supp m$ and fix it. We call a quadruple $(X,d,m,x_0)$ a \emph{pointed metric measure space}. We say that a pointed metric measure space $(X,d,m,x_0)$ is \emph{normalized} if 
 \begin{align}
  \int_{B_1(x_0)}1-d(x_0,\cdot)\,dm=1.\label{eq:normal}
 \end{align}
  For $r\in(0,1)$, define $d_r:=d/r$ and 
 \begin{align}
  m^x_r:=\left(\int_{B_r(x)}1-d_r(x,\cdot)\,dm\right)^{-1}m.\notag
 \end{align}
 Note that the pointed metric measure space $(X,d_r,m^x_r,x)$ is normalized. 
 
 Let $C(\cdot):[0,\infty)\rightarrow [1,\infty)$ be a nondecreasing function. Define $\mathcal{M}_{C(\cdot)}$ the family of pointed metric measure spaces $(X,d,m,\bar x)$ that satisfy 
 \begin{align}
  m(B_{2r}(x))\leq C(R)m(B_r(x))\notag
 \end{align}
 for any $x\in \mathrm{supp}\,m$, and any $0<r\leq R<\infty$. Gigli, Mondino, and Savar\'e have proven that there exists a distance function $\mathcal{D}_{C(\cdot)}:\mathcal{M}_{C(\cdot)}\times \mathcal{M}_{C(\cdot)}\rightarrow [0,\infty]$, which induces the same topology as the Gromov-Hausdorff one on $\mathcal{M}_{C(\cdot)}$ (\cite{GMS}). It is known that every $RCD^*(K,N)$ space for given $K\in\R$, $N\in (1,\infty)$ belongs to $\calM_{C(\cdot)}$ for a common function $C(\cdot):(0,\infty)\rightarrow [1,\infty)$ (see \cite{Stmms2,FS}), more precisely, they satisfy 
 \begin{align}
  \frac{m(B_R(x))}{m(B_r(x))}\leq \frac{\int_0^R\bbS_{K,N}^{N-1}(t)\,dt}{\int_0^r\bbS_{K,N}^{N-1}(t)\,dt}\notag
 \end{align}
 for any $x\in\supp m$, $0<r\leq R$, where 
 \begin{align}
  \bbS_{K,N}(t):=\begin{cases}
                           \sin\left(t\sqrt{\frac{K}{N-1}}\right)&\text{if }K>0,\\
                           t&\text{if }K=0,\\
                           \sinh\left(t\sqrt{\frac{-K}{N-1}}\right)&\text{if }K<0. 
                          \end{cases}\notag
 \end{align}
 Let $(X,d,m)$ be an $RCD^*(K,N)$ space for $K\in\R$ and $N\in(1,\infty)$. By a simple calculation, we have $(X,d_r,m_r^x)$ for some $x\in \supp m$ being an $RCD^*(r^2K,N)$ space. Take a point $x\in\supp m$ and fix it. Consider the family of normalized metric measure spaces $\{(X,d_r,m^x_r,x)\}_{r\in(0,1)}$. The following comes from a generalization of Gromov's compactness theorem. 
 \begin{thm}
  The family of normalized metric measure spaces $\{(X,d_r,m^x_r)\}_{r\in(0,1)}$ is compact with respect to the pointed measured Gromov-Hausdorff topology. Moreover every limit space $(X,d_{r_n},m^x_{r_n},x)\rightarrow (Y,d_Y,m_Y,y)$ is a normalized $RCD^*(0,N)$ space for the non-increasing sequence $\{r_n\}_n$ with $r_n\rightarrow 0$. 
 \end{thm} 
 We define the tangent cone at a point $x\in\supp m$ by 
 \begin{align}
  \mathrm{\mathrm{Tan}}(X,d,m,x):=\left\{(Y,d_Y,m_Y,y)\;;\;(X,d_{r_n},m^x_{r_n},x)\rightarrow(Y,d_Y,m_Y,y)\right\},\notag
 \end{align}
 where $\{r_n\}_n$ is a non-increasing sequence converging to 0. For simplicity, we just denote by $\mathrm{\mathrm{Tan}}(X,x)$ instead of $\mathrm{\mathrm{Tan}}(X,d,m,x)$ if there is no confusion.  
 \begin{defn}\label{def:regularsets}
  Let $(X,d,m)$ be an $RCD^*(K,N)$ space for $K\in\R$ and $N\in(1,\infty)$. We call a point $x\in\supp m$ a $k$-\emph{regular point} if $\mathrm{\mathrm{Tan}}(X,x)=\{(\R^k,d_E,\underline{\mathcal{L}}^k,0)\}$, where $\underline{\mathcal{L}}^k$ is the normalized measure at $0$. We denote the set of $k$-regular points by $\RE_k$. 
 \end{defn}
 Mondino-Naber proved the following \cite{MN}. 
 \begin{thm}[\cite{MN}*{Corollary 1.2}]\label{thm:MN}
  Let $(X,d,m)$ be an $RCD^*(K,N)$ space for $K\in\R$ and $N\in(1,\infty)$. Then 
  \begin{align}
   m\left(X\setminus \bigcup_{1\leq k\leq N}\RE_k\right)=0.\notag
  \end{align}
 \end{thm}
 Note that by Theorem \ref{thm:MN}, we know neither the uniqueness of $\RE_k$ nor $m$ being absolutely continuous to the $k$-dimensional Hausdorff measure $\calH^k$ on $\RE_k$. And on Ricci limit spaces, we know $m(X\setminus \RE_l)=0$ for some $1\leq l\leq N$ (see Theorem \ref{thm:nonricci}, \ref{thm:colricci}).
%
%
\section{Possible volume growth on regular sets}
In this section we prove a behavior of measure of small balls on regular sets under some assumption, which is crucial in the later section. 
\par We say a metric measure space $(Y,r,\nu)$ \emph{nontrivial} if $\di\, (Y,r)>0$ and $\nu(Y)\in(0,\infty]$.  
\begin{defn}
 Let $(X,d,m)$ be a metric measure space. We define the following weakly regular set. 
 \begin{align}
  &\underline{\calW\calE}_k:=\left\{x\in X\;;\;\text{There exists a nontrivial }X'\text{ with }X'\times\R^k\in \mathrm{Tan}(X,x)\right\}.\notag
 \end{align}
\end{defn}
 The author and Lakzian proved the following lemma, which is independently proven by Kell \cite{Kell}. 
 \begin{lem}[\cite{KL}*{Corollary 5.5}, \cite{Kell}*{Lemma 7}]\label{lem:linear}
  Let $(X,d,m)$ be an $RCD^*(K,N)$ space for $K\in\R$, $N\in(1,\infty)$. Assume $\mathrm{supp}\,m=X$ and $X$ is nontrivial. Fix a point $y\in X$. Then for any $R>0$, there exists a constant $C=C(R,y)$ such that 
  \begin{align}
   m(B_s(x))\leq Cs\notag
  \end{align}
  for each $x\in B_R(y)$ and any $s\in(0,1]$. 
 \end{lem}
\begin{prop}
 Let $(X,d,m)$ be an $RCD^*(K,N)$ space for $K\in\R$ and $N\in(1,\infty)$. Take a point $x\in\underline{\calW\calE}_{k-1}$. Suppose that 
 \begin{align}
  \limsup_{r\rightarrow 0}\frac{m(B_r(x))}{r^{\alpha}}=:A<\infty\notag
 \end{align}
 for some $\alpha<k$. Then 
 \begin{align}
  \liminf_{r\rightarrow 0}\frac{m(B_r(x))}{r^{\alpha}}=0.\notag
 \end{align} 
\end{prop}
\begin{proof}
 Let $\{r_i\}_{i\in\mathbb{N}}$ be a sequence that realizes $X_{r_i}\rightarrow X'\times\R^{k-1}\in \mathrm{Tan}(X,x)$, where $(X',d',m')$ is nontrivial. Since $(X',d',m')$ is an $RCD^*(0,N-k+1)$ space by the splitting theorem \cite{Gsplit}, $m'(B_r(x'))\leq Cr$ for small $r>0$ holds by Lemma \ref{lem:linear}. Thus $\calL^{k-1}\times m'(B_r(0,x'))\leq \calL^{k-1}(B_{r}(0))m'(B_{r}(x'))\leq C_0r^k$ holds for some constant $C_0>0$. For any $\epsilon>0$, we obtain 
 \begin{align}
  \liminf_{r\rightarrow 0}\frac{m(B_r(x))}{r^{\alpha}}&\leq \liminf_{i\rightarrow\infty}\frac{m(B_{\epsilon r_i}(x))}{(\epsilon r_i)^{\alpha}}\notag\\
  &=\liminf_{i\rightarrow \infty}\frac{m^x_{r_i}(B^{d_{r_i}}_{\epsilon}(x))}{\epsilon^{\alpha}}\cdot\frac{\int_{B_{r_i}(x)}1-\frac{1}{r_i}d(x,\cdot)\,dm}{r_i^{\alpha}}\notag\\
  &\leq \liminf_{i\rightarrow \infty}\frac{m^x_{r_i}(B^{d_{r_i}}_{\epsilon}(x))}{\epsilon^{\alpha}}\cdot \frac{m(B_{r_i}(x))}{r_i^{\alpha}}\notag\\
  &\leq (A+1)\lim_{i\rightarrow\infty}\frac{m^x_{r_i}(B^{d_{r_i}}_{\epsilon}(x))}{\epsilon^{\alpha}}\notag\\
  &\leq (A+1)\frac{\calL^{k-1}\times m'(B_{\epsilon}(0,x'))}{\epsilon^{\alpha}}\notag\\
  &\leq C_0(A+1)\epsilon^{k-\alpha}.\notag
 \end{align}
 Since $\epsilon>0$ is arbitrary small, we have the conclusion. 
\end{proof}
\begin{prop}\label{prop:sup}
 Let $(X,d,m)$ be an $RCD^*(K,N)$ space for $K\in\R$ and $N\in(1,\infty)$. Take a point $x\in\RE_k$. Suppose that 
 \begin{align}
  \liminf_{r\rightarrow 0}\frac{m(B_r(x))}{r^{\beta}}=:D>0\notag
 \end{align}
 holds for some $\beta>k$. Then we obtain
 \begin{align}
  \limsup_{r\rightarrow0}\frac{m(B_r(x))}{r^{\beta}}=\infty.\notag
 \end{align}
\end{prop}
Before proving Proposition \ref{prop:sup}, we show the following lemma. 
\begin{lem}\label{lem:sup}
 Take $r\in(0,1)$ and define the sequence $\{r_n\}_{n\in\mathbb{N}}$ by $r_n:=r^n$. Assume $D\in(0,\infty)$. 
 Then for sufficiently large $n$, 
 \begin{align}
  \int_{B_{r_n}(x)}1-\frac{1}{r_n}d(x,\cdot)\,dm\geq \frac{D}{2}\frac{1-r}{r}r_n^{\beta}\frac{r^{\beta+1}}{1-r^{\beta+1}}\notag
 \end{align}
 holds. 
\end{lem}
\begin{proof}
 This lemma is just a calculation. For sufficiently large $n$, we have $m(B_{r_n}(x))/r_n^{\beta}\geq D/2$. For simplicity, we denote $B_r(x)$ by $B_r$ and $d(x,\cdot)$ by $d_x$. Then 
 \begin{align}
  &\int_{B_{r_n}}1-\frac{1}{r_n}d_x\,dm=\int_{B_{r_n}\setminus B_{r_{n+1}}}1-\frac{1}{r_n}d_x\,dm+\int_{B_{r_{n+1}}}1-\frac{1}{r_n}d_x\,dm\notag\\
  &\geq \int_{B_{r_{n+1}}\setminus B_{r_{n+2}}}1-\frac{1}{r_n}d_x\,dm+\int_{B_{r_{n+2}}}1-\frac{1}{r_n}d_x\,dm\notag\\
  &\geq \left(1-\frac{r_{n+1}}{r_n}\right)m(B_{r_{n+1}}\setminus B_{r_{n+2}})+\int_{B_{r_{n+2}}\setminus B_{r_{n+3}}}1-\frac{1}{r_n}d_x\,dm+\int_{B_{r_{n+3}}}1-\frac{1}{r_n}d_x\,dm\notag\\
  &\geq (1-r)\sum_{i=1}^{\infty}r^{i-1}m(B_{r_{n+i}})\notag\\
  &\geq (1-r)\frac{D}{2}\sum_{i=1}^{\infty}r^{i-1}r_{n+i}^{\beta}\notag\\
  &=\frac{D}{2}\frac{1-r}{r}r_n^{\beta}\sum_{i=1}^{\infty}r^{(\beta+1)i}\notag\\
  &=\frac{D}{2}\frac{1-r}{r}r_n^{\beta}\frac{r^{\beta+1}}{1-r^{\beta+1}}.\notag
 \end{align}
\end{proof}
\begin{proof}[Proof of Proposition \ref{prop:sup}]
Nothing to prove if $D=\infty$, so assume $D\in(0,\infty)$. Since $D<\infty$, by Lemma \ref{lem:sup}, 
 \begin{align}
  \int_{B_{r_n}(x)}1-\frac{1}{r_n}d(x,\cdot)\,dm\geq \frac{D}{2}\frac{1-r}{r}r_n^{\beta}\frac{r^{\beta+1}}{1-r^{\beta+1}}\notag
 \end{align} 
 holds. Note that since $x\in\RE_k$, for any $\epsilon>0$, $m^x_{r_n}(B_{\epsilon}^{d_{r_n}}(x))\rightarrow C_1\epsilon^k$ for some constant $C_1$, which is independent of the choice of the sequence $\{r_n\}_n$. We obtain 
 \begin{align}
  \limsup_{r\rightarrow 0}\frac{m(B_r(x))}{r^{\beta}}&\geq \liminf_{n\rightarrow\infty}\frac{m(B_{\epsilon r_n}(x))}{(\epsilon r_n)^{\beta}}\notag\\
  &=\liminf_{n\rightarrow \infty}\frac{m^x_{r_n}(B^{d_{r_n}}_{\epsilon}(x))}{\epsilon^{\beta}}\frac{\int_{B_{r_n}(x)}1-\frac{1}{r_n}d(x,\cdot)\,dm}{r_n^{\beta}}\notag\\
  &\geq \liminf_{n\rightarrow\infty}\frac{m^x_{r_n}(B^{d_{r_n}}_{\epsilon}(x))}{\epsilon^{\beta}}\cdot\frac{D}{2}\frac{1-r}{r}\frac{r^{\beta+1}}{1-r^{\beta+1}}\notag\\
  &=\frac{C_1D}{2}\frac{1-r}{r}\frac{r^{\beta+1}}{1-r^{\beta+1}}\epsilon^{-(\beta-k)}.\notag
 \end{align}
 Since $\epsilon>0$ is arbitrary, we have the conclusion. 
\end{proof}
%
%
\section{Bishop inequality and regular sets on $RCD$ spaces}
From now, we always assume that every metric measure space $(X,d,m)$ satisfies $\supp m=X$. 
In this section we define variants of the Bishop inequalities on $RCD^*(K,N)$ spaces. We prove every $RCD^*(K,N)$ space satisfies \emph{the weak generalized Bishop inequality}. For any $K'\in\R$ and $N'\in \mathbb{N}$, define $\underline{V}_{K',N'}$ as the normalized volume measure on the space form of curvature $K'$ and of dimension $N'$, where normalized means in the sense of (\ref{eq:normal}). 
\begin{defn}\label{def:genBI}
 Let $(X,d,m)$ be a metric measure space. For given $K'\in\R$ and $N'\in\mathbb{N}$, we say that $(X,d,m)$ satisfies \emph{the generalized Bishop inequality} $BI(K',N')$, if for any $\epsilon>0$, and any $x\in X$, there exists a small number $r_{\epsilon,x}>0$ such that 
  \begin{align}
   m_r^x(B^{d_r}_s(x))\leq (1+\epsilon)\underline{V}_{r^2K',N'}(s)\label{eq:genBI}
  \end{align} 
  holds for any $r\in(0,r_{\epsilon,x})$ and any $s\in(0,\min\{\di\, X, 1\})$. 
  \end{defn}
\begin{rem}
 It is obvious that $N$-dimensional Riemannian manifolds with Ricci curvature bounded below by $K$ and non-collapsing $(K,N)$-Ricci limit spaces satisfy the generalized Bishop inequality $BI(K,N)$. However a collapsing Ricci limit space does not satisfy $BI(K,N)$ inequality since at $m$-a.e. points, the behavior of volumes of balls are controlled by a dimension in the sense of Colding-Naber, see \cite{CC1,CC2,CC3,CNholder}. 
\end{rem}
The generalized Bishop inequality might be thought that it is too strong condition for general metric measure spaces. However a bit weaker condition is always satisfied for any finite dimensional $RCD$ spaces. 
\begin{prop}
 Let $(X,d,m)$ be an $RCD^*(K,N)$ space. Then $(X,d,m)$ satisfies the following: for any $\epsilon>0$, $x\in\calR_l$ for $1\leq l\leq N$, and $s\in (0,\min\{\di\, X,1\})$, there exists a small number $r_{\epsilon,x,s}>0$ such that 
 \begin{align}
  m^x_r(B^{d_r}_s(x))\leq (1+\epsilon)\underline{V}_{r^2K,l}(s)\notag
 \end{align}  
 holds for any $r\in(0,r_{\epsilon,x,s})$. 
\end{prop}
\begin{proof}
 Suppose $(X,d,m)$ does not satisfy the above property. Then there exist $\epsilon>0$, $x\in\RE_l$ for some $1\leq l\leq N$, and $s\in(0,1)$ such that 
 \begin{align}
  m^x_{r_n}(B^{d_{r_n}}_s(x))>(1+\epsilon)\underline{V}_{r_n^2K,l}(s)\notag
 \end{align}
 holds for some decreasing sequence $\{r_n\}_{n\in\mathbb{N}}$ with $r_n\rightarrow 0$. Letting $n\rightarrow\infty$ leads 
 \begin{align}
  \underline{\mathcal{L}^l}(B_s(0))\geq(1+\epsilon) \underline{V}_{0,l}(s)=(1+\epsilon)\underline{\mathcal{L}^l}(B_s(0)).\notag
 \end{align} 
 This is a contradiction. 
\end{proof}
The following is the main result in this paper. 
\begin{thm}\label{thm:noncollapsing}
 Let $(X,d,m)$ be a $RCD^*(K,N)$ space for $K\in\R$ and $N\in\mathbb{N}$. Assume $X$ is nontrivial and satisfies the generalized Bishop inequality. Then $\RE_l=\emptyset$ for any $l<N$. In particular $m(X\setminus \RE_N)=0$. 
\end{thm}
\begin{proof}
 Suppose $\RE_l\neq\emptyset$ for some $1\leq l< N$. Take a point $x\in\RE_l$. By the Bishop-Gromov inequality, 
 \begin{align}
  \liminf_{r\rightarrow0}\frac{m(B_r(x))}{r^N}\geq \frac{m(B_1(x))}{V_{K,N}(1)}\limsup_{r\rightarrow 0}\frac{V_{K,N}(r)}{r^N}>0\notag
 \end{align}
 holds. On the other hand, since $(X,d,m)$ satisfies the generalized Bishop inequality, for small $\epsilon>0$, there exists a small number $r_{\epsilon,x}>0$ that satisfies (\ref{eq:genBI}). Take $0<r'<r_{\epsilon,x}$ and fix it. For any $r>0$, define $s:=r/r'$. We have  
 \begin{align}
  \frac{m(B_r(x))}{r^N}&=\frac{m(B_{r's}(x))}{(r's)^N}\notag\\
  &=\frac{m^x_{r'}(B^{d_{r'}}_s(x))}{s^N}\frac{1}{r'^N}\int_{B_{r'}(x)}1-\frac{1}{r'}d(x,\cdot)\,dm\notag\\
  &\leq (1+\epsilon)\frac{\underline{V}_{r'^2K,N}(s)}{s^N}\frac{m(B_{r'}(x))}{r'^N}\notag\\
  &\rightarrow C'(1+\epsilon)<\infty\notag
 \end{align}
 as $s\rightarrow 0$ for some $C':=C'(r',N)>0$. Here $s\rightarrow 0$ is equivalent to $r\rightarrow 0$. Therefore
 \begin{align}
  \limsup_{r\rightarrow 0}\frac{m(B_r(x))}{r^N}<\infty.\notag
 \end{align} 
 This contradicts Proposition \ref{prop:sup}. Hence $\RE_l=\emptyset$ if $l<N$. Hence combining this fact and Theorem \ref{thm:MN} leads the consequence.  
\end{proof}
The following is a direct consequence of the proof of Theorem \ref{thm:noncollapsing}. 
\begin{cor}\label{cor:ulbound}
 Let $(X,d,m)$ be an $RCD^*(K,N)$ space for $K\in\R$, $N\in\bbN$ and assume $(X,d,m)$ satisfies the generalized Bishop inequality. Then 
 \begin{align}
  0<\liminf_{r\rightarrow0}\frac{m(B_r(x))}{r^N}\leq \limsup_{r\rightarrow0}\frac{m(B_r(x))}{r^N}<\infty\label{eq:ulbound}
 \end{align}
 holds for any $x\in X$. In particular $(\ref{eq:ulbound})$ holds for any $x\in\calR_N$.
\end{cor}
 For given two Borel measures $\mu,\nu$ on $X$, we say $\mu$ and $\nu$ are \emph{equivalent} if $\mu$ is absolutely continuous with respect to $\nu$, $\mu\ll\nu$, and also $\nu\ll\mu$. 
\begin{cor}\label{cor:Hdim}
 Let $(X,d,m)$ be the same as Corollary $\ref{cor:ulbound}$. Then $m$ is equivalent to the $N$-dimensional Hausdorff measure $\calH^N$ on $\calR_N$. Moreover the Hausdorff dimension of $(X,d)$ is $N$. 
\end{cor}
\begin{proof}
 We define the set $C_i$, $i=1,2,\ldots$, by 
 \begin{align}
  C_i:=\left\{ x\in X\;;\;\frac{1}{i}< \liminf_{r\rightarrow0}\frac{m(B_r(x))}{r^N}< i,\;\text{for }r\in(0,1)\right\}.\notag
 \end{align} 
 By Proposition \ref{prop:measurable} in Appendix \ref{app:Borel}, we claim that $C_i$ is a Borel set for each $i\in\mathbb{N}$. 
 By Corollary \ref{cor:ulbound}, $\RE_N\subset \cup_{i\in\mathbb{N}}C_i$. By Theorem 2.4.3 in \cite{AT}, $m\vert_{C_i}$ is equivalent to $\mathcal{H}^N$ on $C_i$. Since $m(X\setminus \RE_N)=0$, there exist a number $i\in\mathbb{N}$ and a measurable subset $A\subset C_i$ such that $0<c_i^{-1}m(A)\leq \mathcal{H}^N(A)\leq c_im(A)<\infty$ hold for some constant $c_i>0$. Thus $\mathrm{dim}_{\mathcal{H}}X\geq N$. On the other hand, by Corollary 2.5 in \cite{Stmms2}, we have $\mathrm{dim}_{\mathcal{H}}X\leq N$ if $\supp m=X$. Thus we have the conclusion. 
\end{proof}
%
%
\section{Tangent cones} 
In this section, we consider tangent cones at any points on $RCD^*(K,N)$ spaces with the generalized Bishop inequality $BI(K,N)$. For noncollapsing Ricci limit spaces, all tangent cones are metric cones. Here we prove the same result in our setting. The following result plays a key role in our proof. 
\begin{thm}[\cite{GDP}*{Theorem 1.1}]
 Let $(X,d,m)$ be an $RCD^*(0,N)$ space with $\supp m=X$, $o\in X$, and $R>r>0$ such that 
 \begin{align}
  m(B_R(o))=\left(\frac{R}{r}\right)^Nm(B_r(o)).\notag
 \end{align}
 Then exactly one of the following holds: 
 \begin{enumerate}
  \item $S_{R/2}(o)$ contains only one point. In this case $(X,o)$ is pointed isometric to $([0,\di\, (X,d)],0)$ and via the isometry, $m\vert_{B_R(o)}$ can be seen as $cx^{N-1}dx$ for $c:=Nm(B_R(o))$.
  \item $S_{R/2}(o)$ contains two points. In this case $(X,d)$ is a 1-dimensional Riemannian manifolds possibly with boundary, and there is a bijective local isometry from $B_R(o)$ to $(-R,R)$ sending $o$ to 0 and the measure $m\vert_{B_R(o)}$ to the measure $c\vert x\vert^{N-1}dx$ for $c:=Nm(B_R(o))/2$. Moreover such local isometry is an isometry when restricted to $\overline{B}_{R/2}(o)$.
  \item $S_{R/2}(o)$ contains more than two points. In this case $N\geq 2$ and there exists an $RCD^*(N-2,N-1)$ space $(Z,d_Z,m_Z)$ with $\di\, (Z,d_Z)\leq \pi$ such that the ball $B_R(o)$ is locally isometric to the ball $B_R(o_Y)$ of the cone built over $Z$. Moreover such local isometry is an isometry when restricted to $\overline{B}_{R/2}(o)$. 
 \end{enumerate}
\end{thm}
\begin{thm}\label{thm:metriccone}
 Let $(X,d,m)$ be an $RCD^*(K,N)$ space for $K\in\R$ and $N\in(1,\infty)$. Assume the generalized Bishop inequality $BI(K,N)$. Then tangent cones are metric cones. Moreover, every tangent cone $(Y,d_Y,m_Y)$ is isometric to a cone over an $RCD^*(N-2,N-1)$ space $(Z,d_Z,m_Z)$. 
\end{thm}
\begin{proof}
 Take an arbitrary point $x\in X$ and fix it. By the Corollary \ref{cor:ulbound}, we have 
 \begin{align}
  0<\limsup_{r\rightarrow0}\frac{m(B_r(x))}{r^N}<\infty.\notag
 \end{align}
 And by the Bishop-Gromov inequality, we have the monotonicity of $r\mapsto m(B_r(x))/V_{K,N}(r)$. Hence we obtain the limit 
 \begin{align}
  C_1:=\lim_{r\rightarrow 0}\frac{m(B_r(x))}{r^N}.\notag
 \end{align}
 Take a tangent cone $(Y,d_Y,m_Y,y)\in \mathrm{Tan}(X,x)$. Denote by $\{r_n\}_{n\in\mathbb{N}}$ the decreasing sequence that realized the convergent sequence $\{(X,d_{r_n},m^x_{r_n},x)\}_{n\in\bbN}$ to $Y$. Since 
 \begin{align}
  \frac{r^N}{m(B_r(x))}\leq \frac{r^N}{\int_{B_r(x)}1-d_r(x,\cdot)\,dm}\leq 2^{N+1}\frac{r^N/2^N}{m(B_{r/2}(x))}\notag
 \end{align}
 holds, we are able to take a subsequence $\{r_{n_k}\}_{k\in\bbN}$ such that $r_{n_k}^N/\int_{B_{r_{n_k}}(x)}1-d_{r_{n_k}}(x,\cdot)\,dm$ converges to a positive and finite constant $C_2$. For simplicity, we just write $\{r_k\}$ instead of $\{r_{n_k}\}$. Then we obtain 
 \begin{align}
  m_Y(B_s(y))&=\lim_{n\rightarrow \infty}m^x_{r_n}(B^{d_{r_n}}_s(x))=\lim_{k\rightarrow\infty}m^x_{r_k}(B^{d_{r_k}}_s(x))\notag\\
  &=s^N\lim_{k\rightarrow\infty}\frac{m(B_{r_ks}(x))}{r_k^Ns^N}\frac{r_k^N}{\int_{B_{r_k}(x)}1-d_{r_k}(x,\cdot)\,dm}\notag\\
  &=C_1C_2s^N.\notag
 \end{align}
 Therefore the equality 
 \begin{align}
  \frac{m_Y(B_{s_2}(y))}{m_Y(B_{s_1}(y))}=\left(\frac{s_2}{s_1}\right)^N\notag
 \end{align}
 holds for any $s_1,s_2>0$. Since $N>1$ and $\calR=\calR_N$, $(Y,d_Y,m_Y)$ has no 1-dimensional regular set by Theorem 1.1 in \cite{KL}. Therefore, by Theorem \ref{thm:metriccone}, $(Y,d_Y,m_Y)$ has to be a metric cone built over an $RCD^*(N-2,N-1)$ space $(Z,d_Z,m_Z)$.  
 \end{proof}
 \begin{rem}\label{rem:boundaryTan}
  By Theorem \ref{thm:noncollapsing}, Corollary \ref{cor:ulbound}, Theorem \ref{thm:metriccone}, it seems that $RCD$ spaces with the generalized Bishop inequality are able to be called ``Non-collapsing $RCD$ spaces". However the following example holds. Consider the closed convex domain $U$ in the $N$-dimensional Euclidean space. Since non-collapsing Ricci limit space do not have boundaries(see \cite{CC1}), such a domain $U$ is not a non-collapsing Ricci limit space. However, $U$ with the standard Euclidean metric and the Lebesgue measure satisfies $RCD(0,N)$ condition, moreover, it satisfies the generalized Bishop inequality. Still, as far as the author knows, whether $U$ is a collapsing Ricci limit space or not is unknown.  
 \end{rem}
\appendix
\section{Continuity of the volume of balls on metric measure spaces}\label{app:Borel}
Let $(X,d,m)$ be a geodesic metric measure space with $\supp m=X$. We say that $(X,d,m)$ satisfies the \emph{$BG(K,N)$ condition} for $K\in\R$, $N\in(1,\infty)$ if 
\begin{align}
 \frac{m(B_R(x))}{m(B_r(x))}\leq \frac{\int_0^R\bbS_{K,N}(t)^{N-1}\,dt}{\int_0^r\bbS_{K,N}(t)^{N-1}\,dt}\notag
\end{align}
holds for any $x\in X$, $0<r\leq R$. Note that, as already stated in Section \ref{sec:pre}, all $RCD^*(K,N)$ spaces satisfy the $BG(K,N)$ condition. Define the function $F(r):=\int_0^r\bbS_{K,N}^{N-1}(t)\,dt$. In \cite{K}, we have the following estimate. 
\begin{lem}\label{lem:ballasymp}
 Let $(X,d,m)$ be a geodesic metric measure space satisfying $BG(K,N)$ condition for $K\in\R$, $N\in(1,\infty)$. For given $r>0$ we take two points $x,y\in X$ such that $d(x,y)<r$. Then we obtain 
 \begin{align}
  \frac{m(B_r(x)\setminus B_r(y))}{m(B_r(x))}\leq \frac{F'(r-d(x,y)/2)}{F(r+d(x,y)/2)}d(x,y)+o(d(x,y))\notag
 \end{align}
 as $d(x,y)\rightarrow0$. 
\end{lem}
In the same spirit as \cite{K}, we prove the local Lipschitz continuity of the volume of balls. 
\begin{lem}\label{lem:locLip}
 Let $(X,d,m)$ be an $RCD^*(K,N)$ space for $K\in\R$ and $N\in(1,\infty)$. 
 Take a point $x_0\in X$ and fix it. Let $R>0$ be a positive number. Then there exists a positive number $C_2:=C_2(x_0,R)>0$, such that 
 \begin{align}
  \vert m(B_r(x))-m(B_r(y))\vert\leq C_2d(x,y)\notag
 \end{align} 
 for any $x,y\in B_R(x_0)$ and any $r\in(0,1)$. 
\end{lem}
\begin{proof}
 By Lemma \ref{lem:linear}, there exists $C':=C'(x_0,R)>0$ such that $m(B_r(x))\leq C'r$ for any $x\in B_R(x_0)$ and any $r\in(0,1)$. Take two points $u,v\in B_R(x_0)$ with $d(u,v)<r$. For simplicity, we write $d:=d(u,v)$. Then by Lemma \ref{lem:ballasymp}, 
 \begin{align}
  \vert m(B_r(u))-m(B_r(v))\vert&\leq m(B_r(u)\setminus B_r(v))+m(B_r(v)\setminus B_r(u))\notag\\
  &\leq 2\left\{\frac{F'(r-d/2)}{F(r+d/2)}d+o(d)\right\}\left(m(B_r(u))+m(B_r(v))\right)\notag\\
  &\leq \left\{2C'\left(\frac{rF'(r-d/2)}{F(r+d/2)}+\frac{o(d)}{d}\right)\right\}d.\label{lem:locLipm:eq1}
 \end{align}
 Note that the curly bracket in (\ref{lem:locLipm:eq1}) satisfies 
 \begin{align}
   \lim_{r\rightarrow 0}\lim_{d(x,y)\rightarrow 0}\left\{2C'\left(\frac{rF'(r-d/2)}{F(r+d/2)}+\frac{o(d)}{d}\right)\right\}=:\tilde C<\infty.\notag
 \end{align}
 Hence for sufficiently small $d$, (\ref{lem:locLipm:eq1}) is bounded above. 
 Again take arbitrary distinct points $x,y\in B_R(x_0)$. Let $\gamma:[0,1]\rightarrow X$ be a geodesic connecting them. Divide the interval $[0,1]$ into $[0,1]=\cup[t_i,t_{i+1}]$, $i=1,\ldots,n$, so that $d(\gamma_{t_i},\gamma_{t_{i+1}})<r$. By (\ref{lem:locLipm:eq1}), we have
 \begin{align}
  &\vert m(B_r(x))-m(B_r(y))\vert\leq \sum_{i=1}^n\vert m(B_r(\gamma_{t_i}))-m(B_r(\gamma_{t_{i+1}}))\vert\notag\\
  &\leq 2C'\sum_{i=1}^n\left\{\frac{rF'(r-d(\gamma_{t_i},\gamma_{t_{i+1}})/2)}{F(r+d(\gamma_{t_i},\gamma_{t_{i+1}})/2)}+\frac{o(d(\gamma_{t_i},\gamma_{t_{i+1}}))}{d(\gamma_{t_i},\gamma_{t_{i+1}})}\right\}d(\gamma_{t_i},\gamma_{t_{i+1}}).\notag
 \end{align}
 Since we are able to divide $[0,1]$ so that $d(\gamma_{t_i},\gamma_{t_{i+1}})$ takes arbitrary small value, we obtain 
 \begin{align}
  \vert m(B_r(x))-m(B_r(y))\vert\leq \tilde Cd(x,y).\notag
 \end{align}
 Putting $C_2:=\tilde C$ leads the consequence. 
\end{proof}
\begin{prop}\label{prop:measurable}
 Let $(X,d,m)$ be a geodesic metric measure space that satisfies $BG(K,N)$ condition for $K\in\R$, $N\in(1,\infty)$. 
 For given $a,b\in(0,\infty)$ with $a<b$, define the set $C_{a,b}$ by 
 \begin{align}
  C_{a,b}:=\left\{x\in X\;;\;a\leq \liminf_{r\downarrow0}\frac{m(B_r(x))}{r^N}\leq \limsup_{r\downarrow0}\frac{m(B_r(x))}{r^N}\leq b\right\}.\notag
 \end{align} 
 Then $C_{a,b}$ is Borel measurable. 
\end{prop}
\begin{proof}
 Since 
 \begin{align}
  C_{a,b}=\left\{x\in X\;;\;a\leq \liminf_{r\downarrow0}\frac{m(B_r(x))}{r^N}\right\}\cap\left\{x\in X\;;\;\limsup_{\downarrow0}\frac{m(B_r(x))}{r^N}\leq b\right\},\notag
 \end{align} 
 and since the case for the limit supremum are able to be proven by the same way below, we only prove that 
 \begin{align}
  A_a:=\left\{x\in X\;;\;a\leq \liminf_{r\downarrow0}\frac{m(B_r(x))}{r^N}\right\}\notag
 \end{align} 
 is a Borel set. Let $\{r_i\}_{i\in\bbN}$ be a decreasing sequence that satisfies $r_i/r_{i+1}\rightarrow1$ as $i\rightarrow \infty$, for instance take $r_i=i^{-1}$. Consider a family of continuous functions $f_{r_i}(x):=m(B_{r_i}(x))/r_i^N$, where the continuity of each $f_{r_i}$ for $i\in\bbN$ is guaranteed by Lemma \ref{lem:locLip}. Hence $\liminf_if_{r_i}$ is a Borel function. Take an arbitrary $x\in X$ and fix it. Define $\{l_j\}_{j\in\bbN}$ such that $m(B_{l_j}(x))/l_j^N\rightarrow\liminf_{r\downarrow0}m(B_r(x))/r^N$. For each $j\in\bbN$, there exists $i(j)\in\bbN$ such that $r_{i(j)+1}\leq l_j<r_{i(j)}$. Then 
 \begin{align}
  \left(\frac{r_{i(j)+1}}{r_{i(j)}}\right)^N\cdot\frac{m(B_{r_{i(j)+1}}(x))}{r^N_{i(j)+1}}&\leq \left(\frac{r_{i(j)+1}}{l_j}\right)^N\cdot\frac{m(B_{r_{i(j)+1}}(x))}{r^N_{i(j)+1}}\notag\\
  &\leq \frac{m(B_{l_j}(x))}{l_j^N}\leq \frac{m(B_{r_{i(j)}}(x))}{r_{i(j)}^N}\cdot\left(\frac{r_{i(j)}}{l_j}\right)^N\notag\\
  &\leq \frac{m(B_{r_{i(j)}}(x))}{r_{i(j)}^N}\cdot\left(\frac{r_{i(j)}}{r_{i(j)+1}}\right)^N\notag
 \end{align}
 holds for any $j\in\bbN$. By the assumption of $\{r_i\}_{i\in\bbN}$, we obtain 
 \begin{align}
  \liminf_{i\rightarrow\infty}f_{r_i}(x)\leq \lim_{j\rightarrow\infty}\frac{m(B_{l_j}(x))}{l_j^N}=\liminf_{r\downarrow0}\frac{m(B_r(x))}{r^N}\leq \liminf_{i\rightarrow\infty}f_{r_i}(x).\notag
 \end{align}
 Therefore $A_a$ is Borel measurable. 
\end{proof}

\section*{Acknowledgement}
The author would like to thank Professors Shouhei Honda, Shin-ichi Ohta, Tapio Rajala, and Takumi Yokota for their helpful comments and fruitful discussions. He is partly supported by the Grant-in-Aid for JSPS Fellows, The Ministry of Education, Culture, Sports, Science and Technology, Japan and JSPS KAKENHI Grant Number 15K17541. 

\end{document}